# A very accurate method to approximate discontinuous functions with a finite number of discontinuities

E. Stella, C. L. Ladera and G. Donoso

Departamento de Física,

Universidad Simón Bolívar,

Caracas 1086, Venezuela

*clladera@usb.ve*

**Abstract** A simple and very accurate method to approximate a function with a finite number of discontinuities is presented. This method relies on hyperbolic tangent functions of rational arguments as connecting functions at the discontinuities, each argument being the reciprocal of Newton binomials that depend on the abscissae that define the domain of the discontinuous function and upon the abscissae of discontinuities. Our approximants take the form of linear combinations of such hyperbolic tangent functions with coefficients that are obtained by solving a linear system of inhomogeneous equations whose righthand sides are the partition functions that define the given discontinuous function. These approximants are analytic, and being free from the Gibbs phenomenon certainly converge at the discontinuity points much better than other known approximants to discontinuous functions, typical relative errors being of the order of $10^{-14}$ even when as close as $10^{-12}$ to the discontinuity points. Moreover, they can be readily scaled to larger intervals. Our method is here illustrated with a representative set of discontinuous mathematical physics functions, and by studying the dynamics of an oscillator subjected to a discontinuous force, but it can be applied to important cases of discontinuous functions in physics, mathematics, engineering and physical chemistry.

## I Introduction

A number of methods have been developed over the years to approximate important functions in mathematical physics. The list of these methods includes truncated Fourier and truncated Chebyshev series, Fourier-Jacobi and Padé-Jacobi polynomials, Padé-Jacobi rational functions, Padé-Chebyshev and Padé-Legendre polynomials, as well as fractional and quasi-fractional approximations [1-12]. However, the convergence of partial sums of orthogonal series is adversely affected when functions with singularities are approximated. Usually the loss of convergence arises in the regions where the singularities occur, a problem which has come to be known as the Gibbs phenomenon. This phenomenon manifests in an oscillatory behaviour of the approximant in the vicinity of the discontinuity jumps, and constitutes an obstruction to the reconstruction the discontinuous functions (it is in fact an effect that can be even observed as the intensity ripples of Fresnel diffraction patterns when a coherent plane wave diffracts from straight edges [13]). The convergence failure at the discontinuity points has been thoroughly studied for many of those approximation methods. Attrition procedures or modifications are routinely applied to hamper the loss of convergence, *e.g.* the well-known *Lanczos sigma factor* introduced in Fourier series expansions, but only with partial success i.e. invariably the ringing or oscillatory behaviour at the discontinuities persists. It shall be seen that the method of approximation presented in this work is not affected by this lack of convergence at the discontinuity points. Our method for generating approximants to arbitrary discontinuous functions with a finite number of discontinuities is based on the definition of connecting tangent hyperbolic functions at each discontinuity point of the discontinuous function being considered. Our approximants take the form of linear combinations of such connecting functions, with coefficients that are found after solving a simple linear system of inhomogeneous equations, in which the role of independent terms is taken by the partition functions of the initially given piecewise continuous function, the one that we wish to approximate. The approximants we hereby obtain are very accurate, the relative errors being truly small, even in small neighbourhoods about the discontinuity points, and even at points where the approximated function would show asymptotic behaviour. Another advantage of our method is that a single approximant is obtained, one that accurately approximates and represents all the partition functions of the initial piecewise continuous functions. Better even, our approximants are analytic and can be easily scaled *i.e.* once derived for any interval $(x_0,x_f)$ of the independent variable the approximant is good for any scaled interval $(kx_0,kx_f)$, $k$ being any positive real number.

Our approximation method can be applied to the familiar set of discontinuous functions of mathematical physics, and to the solution of motion differential equations with discontinuous righthand sides, *i.e.* to the solution of the so-called *Filippov´ problems* [14, 17]. Such kind of differential equations represent the dynamics of important physics and engineering systems such as motion under dry friction, oscillations in visco-elastic media, electrical circuits with switches or small inductivities, motion in non-differentiable potentials, brake processes with locking phase and optimal control of uncertain systems, just to mention a few. As a matter of fact the discussion whether there exists or not discontinuous functions in physics has not been settled. An excellent criterion to settle the discussion, applicable to the case of physical systems that can operate in different modes is to compare the time scale of the transition from one mode to another with the scale of the dynamics of the individual modes [17]: If the first is much smaller than the second it may then be very advantageous to model the transition as being instantaneous, and such cases are therefore amenable to be tackled with our approximation method. Cases of the real world where the state variable is the position instead of time, and that can be represented by discontinuous functions are also known, *e.g.* the segmentation problem in computer vision which consists in computing a decomposition on a set of regions $R_i$ of the domain $R = \cup R_i$ of an image $g$ such that: (a) the image $g$ varies smoothly and/or slowly within each $R_i$ (b) the image $g$ varies discontinuously and/or rapidly across most of the boundaries $\Gamma_i$ between consecutives $R_i$´s. This one is a problem amenable to be treated with our approximation method. In a different context, an important and fundamental set of discontinuous functions are the ones that usually arise when representing the behaviour of thermodynamic systems as the latter undergo phase transitions. In such transitions the properties of the system, or medium, often change discontinuously, as a result of the change of some independent external parameter, typically the temperature or the heat capacity [18, 19]. For instance our method can be applied to represent the heat capacity of a substance in its different phases using just a single accurate approximating function of temperature at constant pressure, and with this function at hand to calculate the enthalpy, entropy and free energy of the substance.

## II Theory of the new analytic approximating method to piecewise continuous functions

Consider a piecewise continuous function $\boldsymbol{\Psi} = \cup \psi_j$ partitioned into a set of known continuous functions $\{\psi_j\}$ at a finite number $N$ of known discontinuity points $\{x_j\}$ of its real domain $(x_0, x_f)$ (a complex domain may be also considered). That is, the discontinuous function $\boldsymbol{\Psi}$ is defined so that in each of the intervals between pairs of consecutive discontinuity points $(x_j, x_{j+1})$ it is given by the partition continuous function $\psi_j$, as illustrated in Fig. 1. Thus, for $N=2$ we will be dealing with two discontinuity points $x_1$, $x_2$, three partition intervals $(x_o, x_1)$, $(x_1, x_2)$, $(x_2, x_f)$ and three partition functions $\{\psi_1, \psi_2, \psi_3\}$. The aim of this work is to find an accurate, analytic yet simple approximating function $\boldsymbol{\Omega}$ to the discontinuous function $\Psi$ across its entire domain.

Let us begin by defining what we call the *discontinuity-matching* or *connecting* functions $\chi_{sld}$ at the known discontinuity abscissae $x_j$ mentioned above (plotted in Fig.

1). These connecting functions $\chi_{sld}$ are in fact continuous functions which we have chosen to define at each discontinuity point $x_j$ in terms of the hyperbolic tangent function in the following way:

$$\chi_{SLD}(x, x_0, x_i, x_f) \equiv \tanh\left\{\frac{(x_f - x_i)^5 (x - x_0)}{[(x-x_i)^2 + b^2]\,[(x-x_0)^2 + b^2]\,\left[(x-x_f)^2 + b^2\right]}\right\}, x \in R \quad (1)$$

where $(x_i, x_f)$ is the domain of the real variable $x$, $x_0$ is the exact point of discontinuity, when the parameter $b$ go to zero. In this sense $b$ is a control parameter related to the form that the discontinuity is approximated. This parameter may be written as:

$$b = (x_f - x_i).\,k,$$

where $k$ is a constant that is chosen according with the computational power of the computer used. A typical value of $k$ would be $10^{-5}$ for a 64 bytes computer.

Let us now consider a piecewise continuous function $\Psi = U\psi_j$, as the one plotted in Fig. 1.

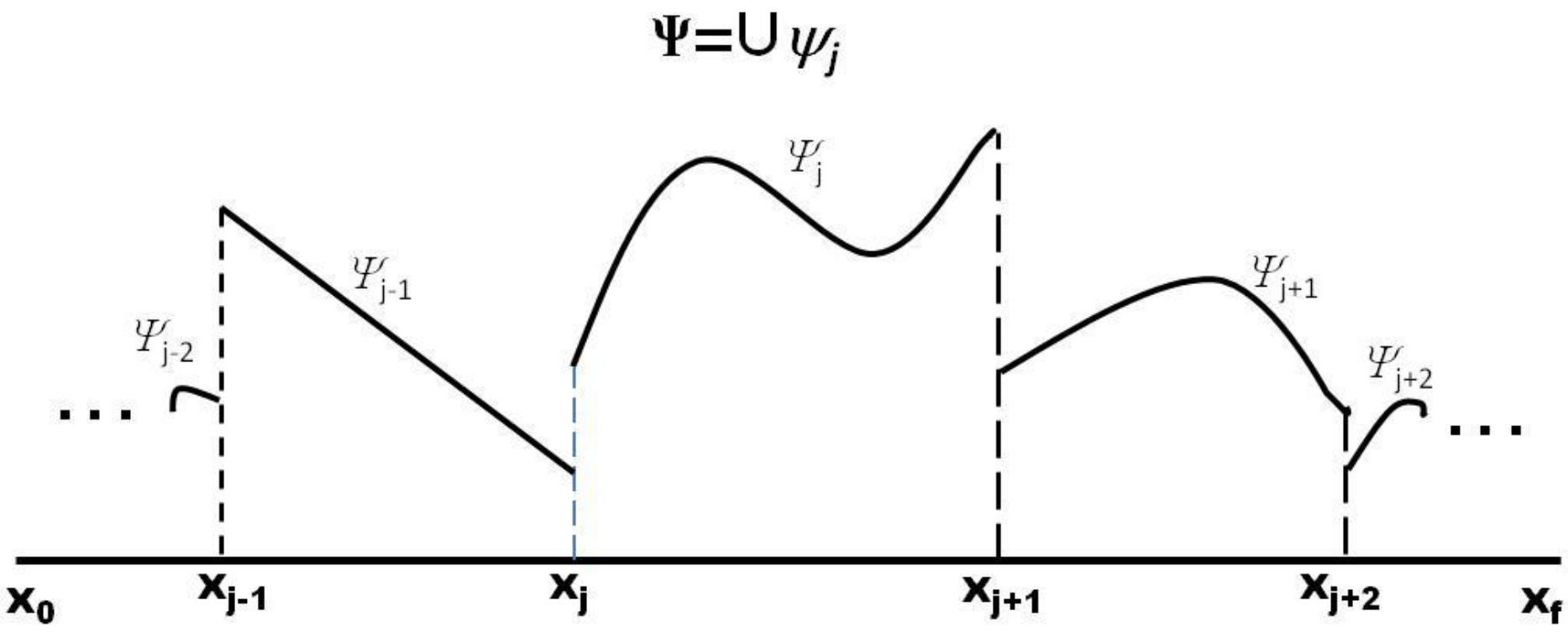

**Fig.1** A piecewise continuous function $\mathbf{\Psi}$ with a finite number $N$ of known discontinuities at the points $\{x_j\}$ is defined on a given real domain $(x_0, x_f)$ by a set of known but arbitrary continuous partition functions $\{\psi_j\}$ defined in their respective partition intervals $(x_j, x_{j+1})$ so that $\mathbf{\Psi} = \cup\psi_j$.

where $(x_f - x_0)$ is the size of the entire domain of the discontinuous function to be approximated. In Appendix 2 we show that our connecting function $\chi_{sld}$ is derivable, and plot the its first derivative (fig. A-2).

It is clear that the number of connecting functions $\chi_{sld}$ is one less than the number of partitions of the given discontinuous function $\mathbf{\Psi} = \cup\psi_j$ (see Fig. 1). The reader may have also noticed that the actual form of the argument of the hyperbolic function in Eq. (1) has been chosen to be rational. It resembles the fractional or quasi-fractional approximants based on truncated series that have been frequently used to approximate continuous functions ([10-12]. The alert reader may have also noted that as defined by Eq. (1) the connecting function is not sufficiently regular for all purposes. A regularized form of our connecting function we present in Appendix A.

As already mentioned in the Introduction we expect our approximant $\boldsymbol{\Omega}$ to the given piecewise continuous function $\boldsymbol{\Psi}$ to be written in the form of a linear combination of the connecting hyperbolic tangent functions $\chi_j$ (see Eq. (1)) with coefficients to be found as explained below in this section. To that effect let us begin by expressing what would eventually become our *approximating function* $\boldsymbol{\Omega}$ in terms of another function

$$\boldsymbol{\Omega}(x) = \int_{x_0}^{x_f} \delta(y-x)\Phi(x,y)dy, \tag{2}$$

a function $\Phi$ which we here introduce and conveniently define in terms of a set of new *auxiliary functions* $\{F_n\}$ of the independent variable $x$, by writing (with $x_{N+1} \equiv x_f$):

$$\begin{cases} \Phi(x,y) \equiv F_0(x) + \sum_{j=1}^{N} F_j(x)\chi_j\big(y,x_j\big), \text{if } x \wedge y \in [x_0,x_{N+1}], & (3a) \\ \Phi(x,y) \equiv 0, \quad \text{otherwise} & (3b) \end{cases}$$

These auxiliary functions $F_n$ are not arbitrary ones; instead they are to be readily found in a self-consistent way as shall be explained a few steps below. We now make the key *ansatz* of our approximation method: the local average of the function $\Phi(x,y)$ over the interval $(x_n,x_{n+1})$ is to be given by the following equation

$$\psi_n(x) = \frac{\int_{x_n}^{x_{n+1}} \Phi(x,y)dy}{x_{n+1}-x_n}\ , \tag{4}$$

i.e. such functional average of the function $\Phi$ must coincide with the known partition continuous function $\psi_j$ that corresponds to the partition interval being considered in the denominator of this equation. With this *ansatz* we may write the integral of Eq. (3a) as,

$$\frac{\int_{x_n}^{x_{n+1}} F_0(x)dy}{x_{n+1}-x_j} + \sum_{j=1}^{N} F_j(x)\frac{\int_{x_n}^{x_{n+1}} \chi_j(y,x_j)\,dy}{x_{n+1}-x_n} = \psi_n(x), \tag{5}$$

or just simply as

$$F_0(x)\frac{\int_{x_n}^{x_{n+1}} dy}{x_{n+1}-x_j} + \sum_{j=1}^{N} F_j(x)\frac{\int_{x_n}^{x_{n+1}} \chi_j(y,x_j)\, dy}{x_{n+1}-x_n} = \psi_n(x), \tag{6}$$

an expression that will allow us to obtain each of the auxiliary functional coefficients $F_j$ introduced above in Eq. (3).

The left-hand side of Eq. (6) may be further simplified by noting, in its two left hand terms, the presence of the averages of the connecting functions $\chi_j$ in the intervals $(x_n, x_{n+1})$ where the pertinent partition function $\psi_n$ is defined. These averages are of course nothing but numbers that we may denote as:

$$S_{n,j} = \frac{\int_{x_n}^{x_{n+1}} \chi_j(y, x_j) dy}{x_{n+1}-x_n} \quad if\ \ j \neq 0 \tag{7a}$$

$$S_{n,0} = \frac{\int_{x_n}^{x_{n+1}} dy}{x_{n+1}-x_n} = 1, \qquad \forall n \tag{7b}$$

Thus, the coefficient $S_{1,3}$ represents the average of the third connecting function $\chi_3$ in the first partition interval $(x_0, x_1)$, and $S_{2,1}$ represents the average of the first connecting function $\chi_1$ in the second partition interval $(x_1, x_2)$, and so on. With these definitions of the numbers $S_{n,j}$ our previous Eq. (6) takes the simple and compact form

$$S_{n0}(x)F_0(x) + \sum_{j=1}^{N} S_{n,j}\, F_n(x) = \psi_n(x), \quad j = 0, 1, 2, \dots N \tag{8}$$

an equation which closely examined represents a linear system of $(N+1)$ inhomogeneous equations in the yet unknown auxiliary functions $\{F_0, F_1, \dots F_N\}$ we just introduced above, the role of independent terms being taken by the partition functions $\psi_j$ that define the given discontinuous function $\boldsymbol{\Psi}$. By solving the linear system of equations (8) the auxiliary functions $\{F_j\}$ are found, as promised above, and our sought approximant $\boldsymbol{\Omega}$ to the given discontinuous function can be finally written using Eqs. (2):

$$\Omega(x) = F_0(x) + \sum_{j=1}^{N} F_j(x)\chi_j(x_j)\,, \tag{9}$$

this being the final step of our method.

As a way to illustrate how simple it is to apply the method just described above let us consider a discontinuous function partitioned into a triplet of known continuous functions $\{\psi_0, \psi_1, \psi_2\}$ and separated by $N = 2$ discontinuity points. Once we define the two required connecting functions $\chi_1$, $\chi_2$ at the known discontinuity abscissae $x_1, x_2$, by using Eq. (1), we can then readily evaluate the nine required coefficients $S_{nj}$ by simply averaging them after Eqs. (7a) and (7b). In the present case these coefficients come up to be

$$\begin{pmatrix} S_{00} & S_{01} & S_{02} \\ S_{10} & S_{11} & S_{12} \\ S_{20} & S_{21} & S_{22} \end{pmatrix} = \begin{pmatrix} 1 & -1 & -1 \\ 1 & 1 & -1 \\ 1 & 1 & 1 \end{pmatrix} \tag{10}$$

And now we may write the linear system of inhomogeneous equations in the three required unknown auxiliary functions $\{F_0, F_1, F_2\}$:

$$\begin{pmatrix} S_{00} & S_{01} & S_{02} \\ S_{10} & S_{11} & S_{12} \\ S_{20} & S_{21} & S_{22} \end{pmatrix} \begin{pmatrix} F_0(x) \\ F_1(x) \\ F_2(x) \end{pmatrix} = \begin{pmatrix} 1 & -1 & -1 \\ 1 & 1 & -1 \\ 1 & 1 & 1 \end{pmatrix} \begin{pmatrix} F_0(x) \\ F_1(x) \\ F_2(x) \end{pmatrix} = \begin{pmatrix} \psi_0(x) \\ \psi_1(x) \\ \psi_2(x) \end{pmatrix}. \tag{9}$$

Once this system (9) of equations is solved in terms of the initially given partition functions $\{\psi_0, \psi_1, \psi_2\}$, the auxiliary functions $\{F_0, F_1, F_2\}$ now being known, the sought approximant is finally written as:

$$\mathbf{\Omega}(x) = F_0(x) + F_1(x)\chi_1(x) + F_2(x)\chi_2(x). \tag{11}$$

The following section is devoted to illustrate how our method of approximation is applied to a particular set of well-known discontinuous functions that recurrently arise in mathematics, physics and engineering, plus to an example of a rather arbitrarily defined discontinuous function. This section on applications we will finish by obtaining an approximant for an interesting mechanical case, taken from Landau and Lifschitz theoretical physics books [20]; it is in fact an example of discontinuous dynamics in which the righthand side of the motion equation is a discontinuous function *i.e.* a Filippov´s problem [14, 15]. In Appendix A we shall present a better and more rigorous definition, in fact what we call a *regularized version* of the key connecting functions of our method. We do *regularize* our connecting functions in order for the approximants to be well defined even at the discontinuity points, and therefore enabling their direct replacement in whatever computational task and when solving differential equations with discontinuous terms [21-23].

## III Applications

This is the section devoted to show how our method of approximation is to be applied. This is firstly done with two well-known mathematical examples and then by way of a discontinuous function which we have arbitrarily defined, and finally with the already mentioned well-known case in mechanics. The quality of the approximants obtained is then assessed by calculating their relative errors, and by closely examining the convergence of the approximant in truly small neighbourhoods at the discontinuity points. The following calculations, in this section, were done with the parameter value $b=0$.

### Example 1 Approximating the *Sign* function

The Sign function has been used for decades as a standard case of discontinuous function, i.e. as a kind of discontinuous function of merit, to evaluate the quality of approximation methods for such kind of functions. Although being a simple case we revisit it here to show the high quality of our approximation method, and to illustrate how simple it is to apply it. Thus consider the real function:

$$S(x) = sign\,(x - 375), \quad x\epsilon(-1000, 1000), \tag{12}$$

plotted in Fig. 2, defined in the domain from $x_0 = -1000$ to $x_f = 1000$, and with a single discontinuity step at the point $x_1 = 375$.

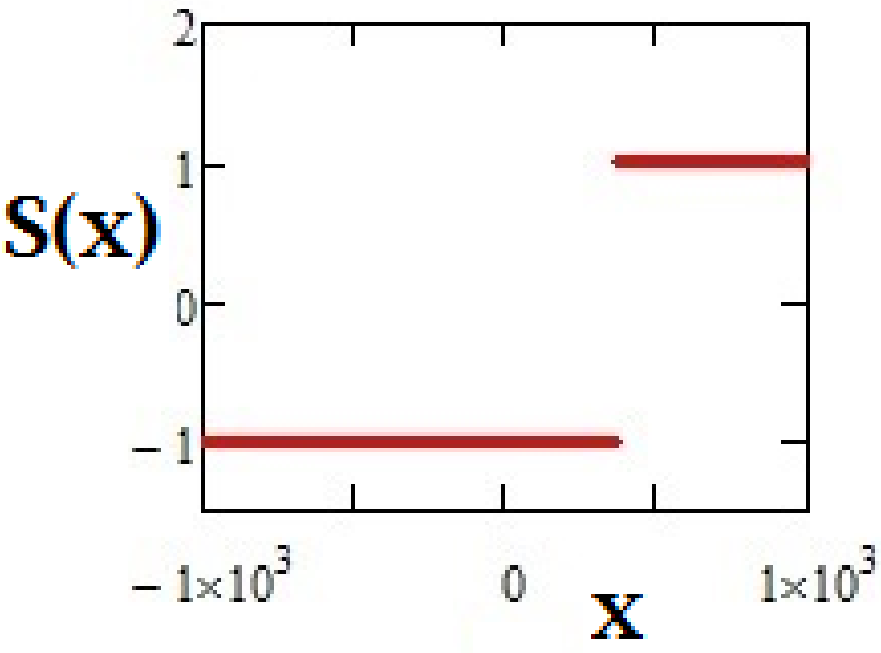

**Fig. 2** Plot of the *Sign* function $S(x) = sign(x\text{-}375)$ whose single discontinuity is at $x = 375$.

This is a piecewise continuous function with only two partition intervals and two partition functions in its domain *(-1000, 1000)*. Using Eq.(1) to define the single connecting function and straightforward application of Eqs.(7), we get the four $S_{nj}$ constants of our approximation method, namely $S_{00} = S_{10} = 1$, plus $S_{01}$, $S_{11}$ given by

$$S_{01} = \frac{1}{1375}\int_{-1000}^{375} tanh\left[\left[\left(\frac{y-x_0}{x_f-x_0}\right)^2\left(\frac{y-x_1}{x_f-x_0}\right)\left(\frac{x_f-y}{x_f-x_0}\right)^2\right]^{-1}\right]dy = -1, \tag{13a}$$

$$S_{11} = \frac{1}{625}\int_{375}^{1000} tanh\left[\left[\left(\frac{y-x_0}{x_f-x_0}\right)^2\left(\frac{y-x_1}{x_f-x_0}\right)\left(\frac{x_f-y}{x_f-x_0}\right)^2\right]^{-1}\right]dy = 1. \tag{13b}$$

With these four constants and the two constant partition functions: $\psi_0 = -1$, $\psi_1 = 1$ (see Fig. (2)) we may write the corresponding linear system of two non-homogeneous equations that will allow us to find the auxiliary functions $F_0$ and $F_1$ (see Eq. (9)). In the current case this system of equations is:

$$\begin{pmatrix} 1 & -1 \\ 1 & 1 \end{pmatrix}\begin{pmatrix} F_0 \\ F_1 \end{pmatrix} = \begin{pmatrix} -1 \\ 1 \end{pmatrix}, \tag{14}$$

whose solutions are $F_0(x) = 0$ and $F_1(x) = 1$. With these two functions and using Eq. (9), the sought approximating function $\boldsymbol{\Omega}$ for the Sign function may be finally written as:

$$\boldsymbol{\Omega}(x) = F_0(x) + F_1(x)\, tanh\left[\left(\frac{x+1000}{2000}\right)^{-2}\left(\frac{x-375}{2000}\right)^{-1}\left(\frac{1000-x}{2000}\right)^{-2}\right],$$

or

$$\boldsymbol{\Omega}(\boldsymbol{x}) = tanh\left[\left(\frac{x+1000}{2000}\right)^{-2}\left(\frac{x-375}{2000}\right)^{-1}\left(\frac{1000-x}{2000}\right)^{-2}\right], \tag{15}$$

which appears plotted in Fig. (3); note the excellent agreement with the discontinuous *Sign* function plotted in Fig. (2).

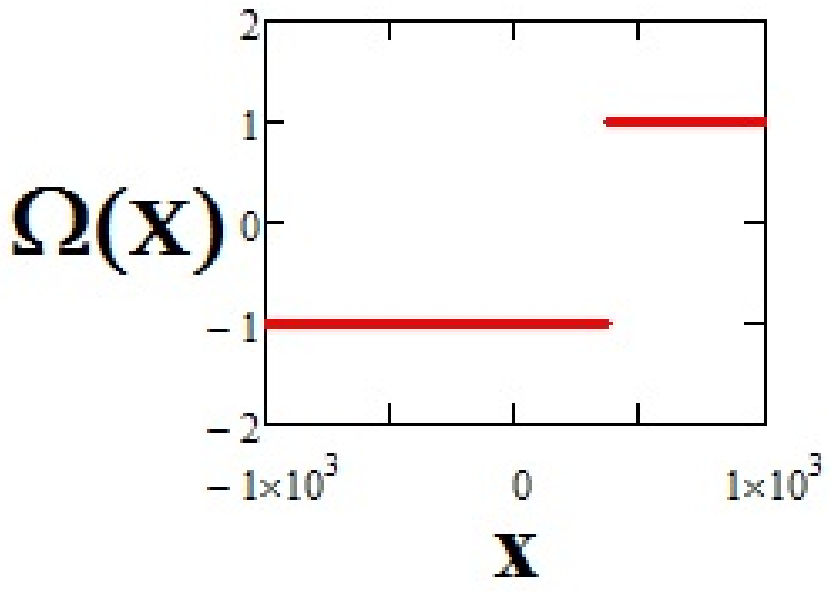

**Fig. 3** Plot of the approximating function *Ω(x)* to the Sign function *S(x) = sign(x-375)* plotted in Fig. 2.

As a matter of fact, and as shown below in Fig (4), the relative error $e(x) = [S(x)-\Omega(x)]/S(x)$ of our approximant to the Sign function is much less than $1 \times 10^{-15}$ both, in its whole domain of definition (Fig 4(a)) and in small neighbourhoods of the single discontinuity point (Fig. 4(b)).

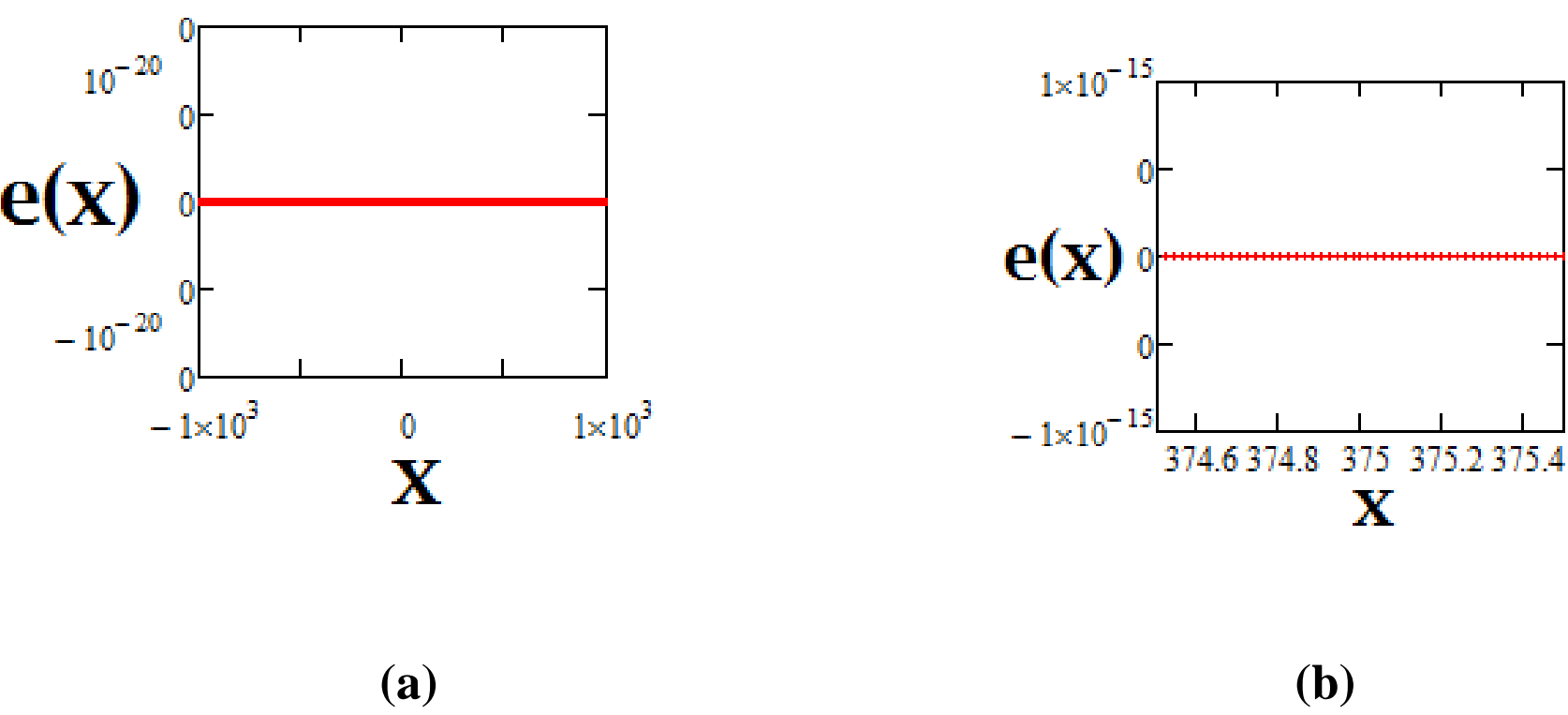

**Fig. 4** Relative error $e(x) = [S(x)-\Omega(x)]/S(x)$ of our approximating function $\Omega$ to the Sign function $S$: (a) in the initial domain $(-1000, 1000)$, (b) about the discontinuity point $x = 375$.

### Example 2. Approximating a *Sawtooth* function

The Sawtooth function is a second example of a simple piecewise continuous function frequently found in mathematics and engineering, and therefore it will serve well to illustrate our approximation method. Thus let us consider the three teeth discontinuous real function $T$ plotted in Fig. 5, that has been defined in the three partition intervals [0,1), [1,2), [2,3) of its domain $[x_0 = 0, x_f = 3]$ in terms of the well-known *floor function*: $T(x) = x - floor(x)$.

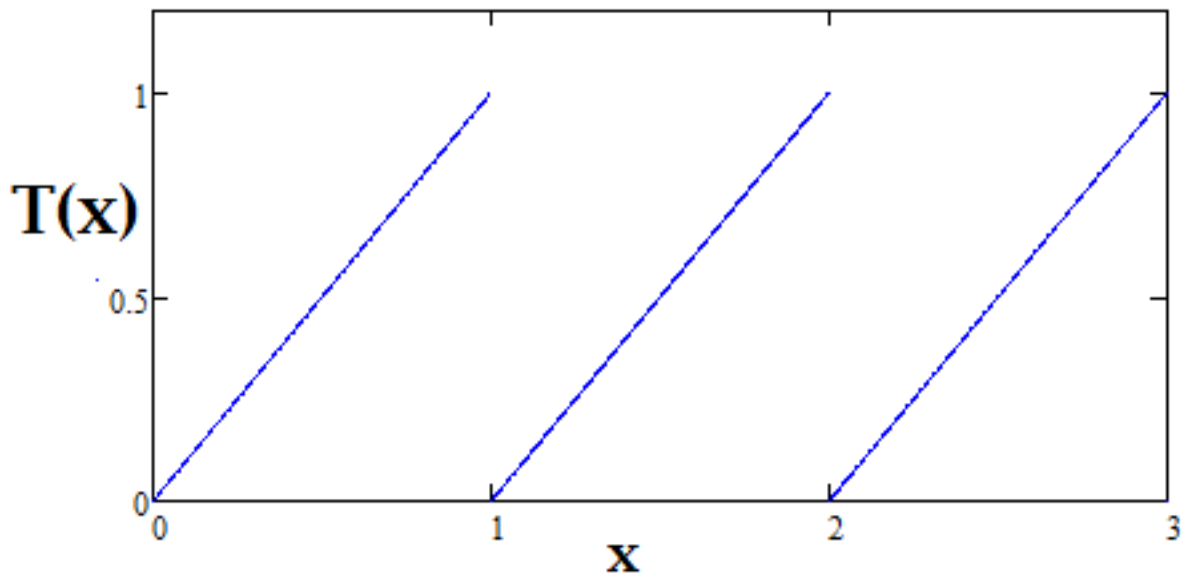

**Fig. 5** Plot of the sawtooth-like function $T(x) = x - floor(x)$ defined in the interval $[x_0=0,\ x_f=3]$ and discontinuous at the abscissae $x_1 = 1$ and $x_2 = 2$.

As shown in the figure this *Sawtooth* is a function defined over three partition intervals, and by three partition functions: $\psi_0(x) = x$ for $x \in [0,1)$, $\psi_1(x) = x-1$ for $x \in [1,2)$, $\psi_2(x) = x-2$ for $x \in [2,3)$. It is discontinuous at the two points $x_1 = 1$ and at $x_2 = 2$.

After straightforward application of Eqs. (1) and (7) of Section 2 we get the nine $S_{nj}$ constant averages of our approximation method, namely $S_{00} = S_{10} = S_{20} = 1$, and $S_{01}$, $S_{11}$, $S_{02}$, $S_{12}$, $S_{21}$ and $S_{22}$ (calculated as done in Example 1 above):

$$\begin{pmatrix} S_{00} & S_{01} & S_{02} \\ S_{10} & S_{11} & S_{12} \\ S_{20} & S_{21} & S_{22} \end{pmatrix} = \begin{pmatrix} 1 & -1 & -1 \\ 1 & 1 & -1 \\ 1 & 1 & 1 \end{pmatrix} \tag{16}$$

With these nine constants and the three given partition functions: $\psi_0$, $\psi_1$, and $\psi_2$ (see Fig. (15)) we may write the corresponding linear system of two inhomogeneous equations to get the auxiliary functions $F_0$ and $F_1$ (see Eq. (9)):

$$\begin{pmatrix} 1 & -1 & -1 \\ 1 & 1 & -1 \\ 1 & 1 & 1 \end{pmatrix} \begin{pmatrix} F_0 \\ F_1 \\ F_2 \end{pmatrix} = \begin{pmatrix} x \\ x-1 \\ x-2 \end{pmatrix}, \tag{17}$$

whose solutions are $F_0(x) = x-1$, $F_1(x) = -1/2$, $F_2(x) = -1/2$. Thus, finally our approximating function Ω(x) for our *Three teeth-Sawtooth* function is given by (see Eq. (9)):

$$\boldsymbol{\Omega}(\mathrm{x}) = (x-1) - \frac{1}{2}\, tanh\left[\left[\left(\frac{x}{3}\right)^2 \left(\frac{x-1}{3}\right)\left(\frac{3-x}{3}\right)^2\right]^{-1}\right] - \frac{1}{2}\, tanh\left[\left[\left(\frac{x}{3}\right)^2 \left(\frac{x-2}{3}\right)\left(\frac{3-x}{3}\right)^2\right]^{-1}\right] \tag{18}$$

which appears plotted in Fig. 6 below.

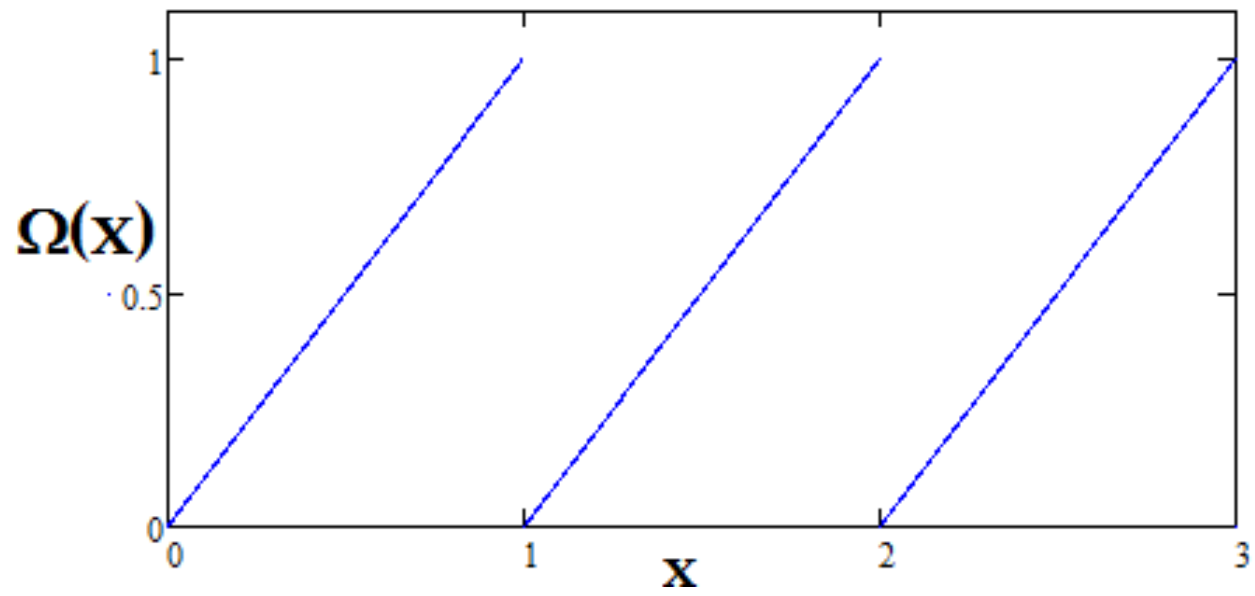

**Fig. 6** Plot of our approximant to the Sawtooth function *T(x) = x- floor(x)* plotted in Fig. 5.

Comparisons of Figs. 5 and 6, show that our approximant reproduces with great accuracy the given sawtooth function *T(x)*. In fact the relative error of our approximant $e(x) = [T(x)-\Omega(x)]/T(x)$ is truly negligible, less than $10^{-13}$ in the neighbourhood of the discontinuity points $x_1 = 1$, $x_1 = 2$, and it is shown in Fig. 7 below.

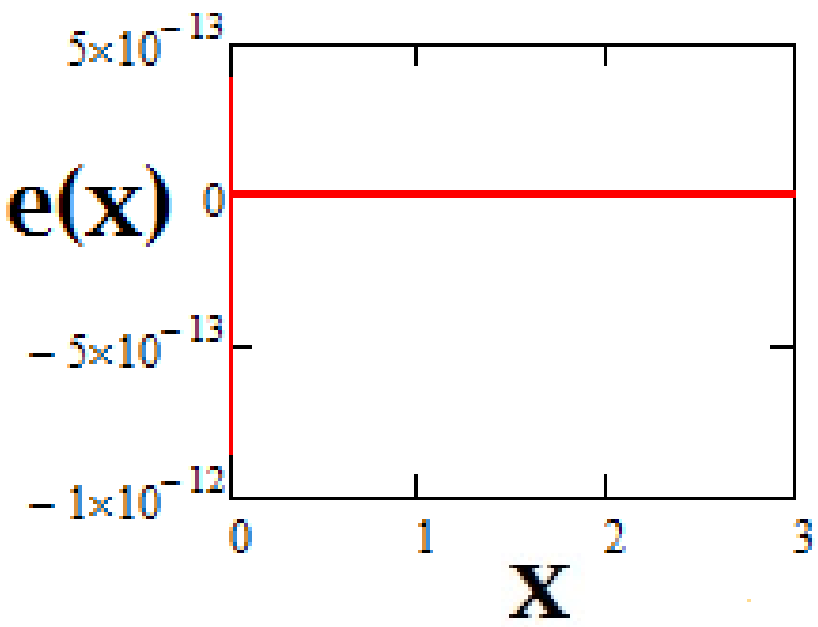

**Fig. 7** Plot of the relative error $e(x)= [T(x)-\Omega(x)]/T(x)$ of our approximant to the *sawtooth* function *T(x) = x- floor(x)* plotted in Fig. 5. It is truly negligible, less than $10^{-13}$ in the whole domain of the function

### Example 3. Approximating an arbitrarily defined discontinuous function

In this example we aim to get an accurate approximant for a piecewise continuous function $\Psi = \cup \psi_j$, $j=1, 2, 3$ defined over the real domain ($x_o = -1.0$, $x_f = 1.0$), discontinuous at the points $x_1 = -0.3$, $x_2 = 0.6$, and arbitrarily defined by the three continuous partition functions given by the vector,

$$\begin{pmatrix} \psi_1(x) \\ \psi_2(x) \\ \psi_3(x) \end{pmatrix} = \begin{pmatrix} ln(1+x^2) \\ \psi_1(x_1)\left(20+20x+5\exp(-4x)\ \sin\left(6\pi\frac{x}{x_2}\right)\right) \\ \psi_2(x_2)+\frac{0.2}{x-0.5}-0.5 \end{pmatrix}, \qquad (19)$$

where $\psi_1(x_1)=0.086177$ and $\psi_2(x_2)=2.770315$, and whose domains are (-1,-0.3), (-0.3, 0.6) and (0.6, 1), respectively. This piecewise continuous function $\mathbf{\Psi}$ appears plotted in Fig. 8 below.

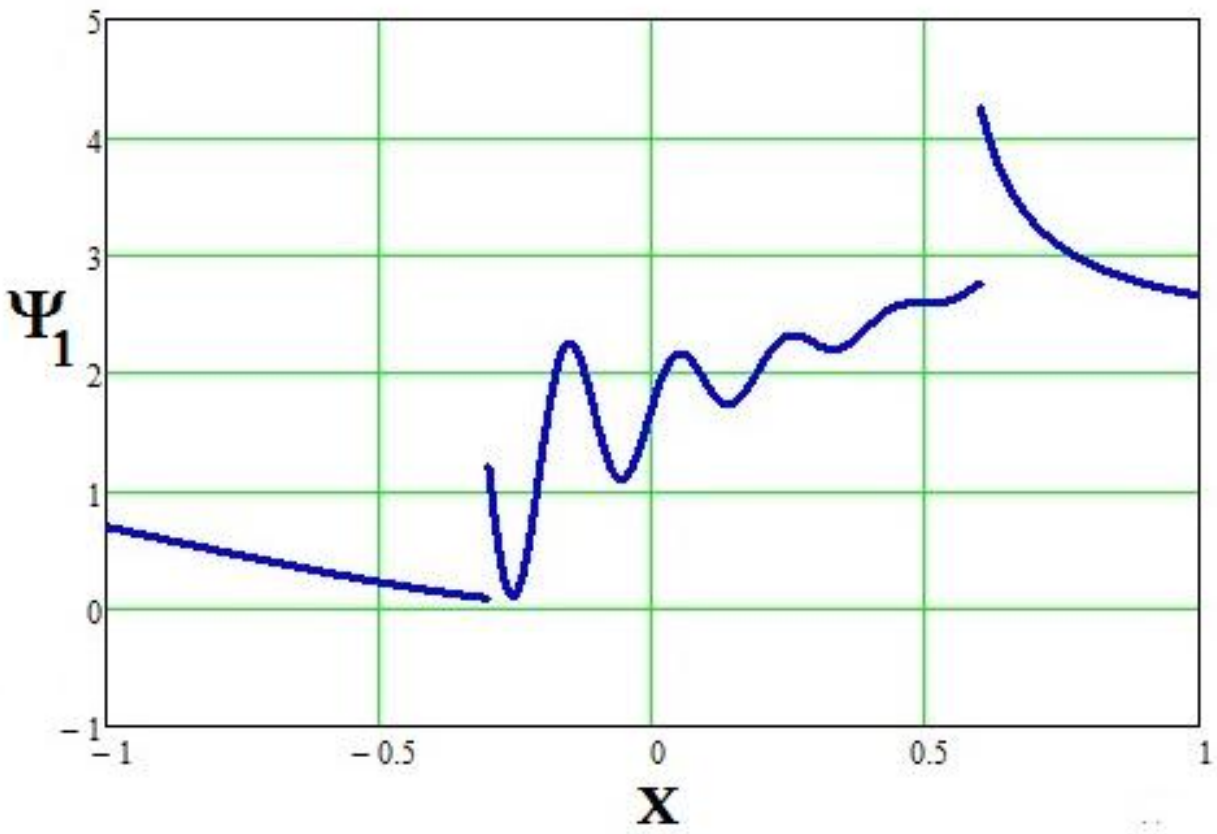

**Fig. 8** Discontinuous function $\mathbf{\Psi} = \cup\psi_j$, $j=1, 2, 3$ to be approximated: it is in fact a piecewise continuous function defined in the real interval [-1,1] with discontinuities at $x_1 = -0.3$ and $x_2 = 0.6$ and partitions functions given in Eq. 19.

The reader may note that the third partition function $\psi_3(x) = \psi_2(x_2) + \frac{0.2}{x-0.5} - 0.5$ in Eq. (19) has been purposefully defined to represent a hyperbola with a vertical asymptote at the point $x=0.5$ of its domain *(0.3, 0.6)*. It shall be seen below that our approximation method does not fail in faithfully representing this discontinuous function in spite of the "pathological" divergence to infinity of the third partition function at $x=0.5$.

As explained in Section 2 one must again begin using Eq. (1) to define the two connecting functions at the two discontinuity points of this example, namely at $x_1=-0.3$, $x_2=0.6$ (recall the domain $x_o= - 1.0$ to $x_f = 1.0$) of the given discontinuous function. For instance at the discontinuity point $x_1=-0.3$ the connecting function is

$$\chi_1(y) \equiv tanh\left[\left[\left(\frac{y+1}{2}\right)^2\left(\frac{y+0.3}{2}\right)\left(\frac{1-y}{2}\right)^2\right]^{-1}\right] \tag{20}$$

The next step in applying our approximation method is the evaluation of the nine coefficients $S_{ij}$, $(i=1,\ 2;\ j=0,\ 1,\ 2)$ over the three partition domains shown in Fig. 2. These we have once again evaluated using Eqs. (7), the results being given in the following matrix,

$$\begin{pmatrix} S_{00} & S_{01} & S_{02} \\ S_{10} & S_{11} & S_{12} \\ S_{20} & S_{21} & S_{22} \end{pmatrix} = \begin{pmatrix} 1 & -1 & -1 \\ 1 & 1 & -1 \\ 1 & 1 & 1 \end{pmatrix}. \tag{21}$$

With these at hand we may now write the linear system of inhomogeneous equations in the three auxiliary functions $F_0$, $F_1$, $F_2$,

$$\begin{pmatrix} 1 & -1 & -1 \\ 1 & 1 & -1 \\ 1 & 1 & 1 \end{pmatrix} \begin{pmatrix} F_0(x) \\ F_1(x) \\ F_2(x) \end{pmatrix} = \begin{pmatrix} \psi_0(x) \\ \psi_1(x) \\ \psi_2(x) \end{pmatrix}, \tag{22}$$

whose solution for the current three partition functions in Eq. (19 ) is as follows:

$$\begin{pmatrix} F_0(x) \\ F_1(x) \\ F_2(x) \end{pmatrix} = \begin{pmatrix} 1/2\,[\psi_0(x) + \psi_2(x)] \\ 1/2\,[-\psi_0(x) + \psi_1(x)] \\ 1/2\,[-\psi_1(x) + \psi_2(x)] \end{pmatrix}, \tag{23}$$

And with these three known auxiliary functions $F_0$, $F_1$ and $F_2$ the sought approximant $\boldsymbol{\Omega_1}$ may be finally written using Eq. (10) as

$$\Omega_1(x) = F_0(x) + F_1(x)\chi_1(x) + F_2(x)\chi_2(x), \tag{24}$$

which appears plotted in Fig. 9 below.

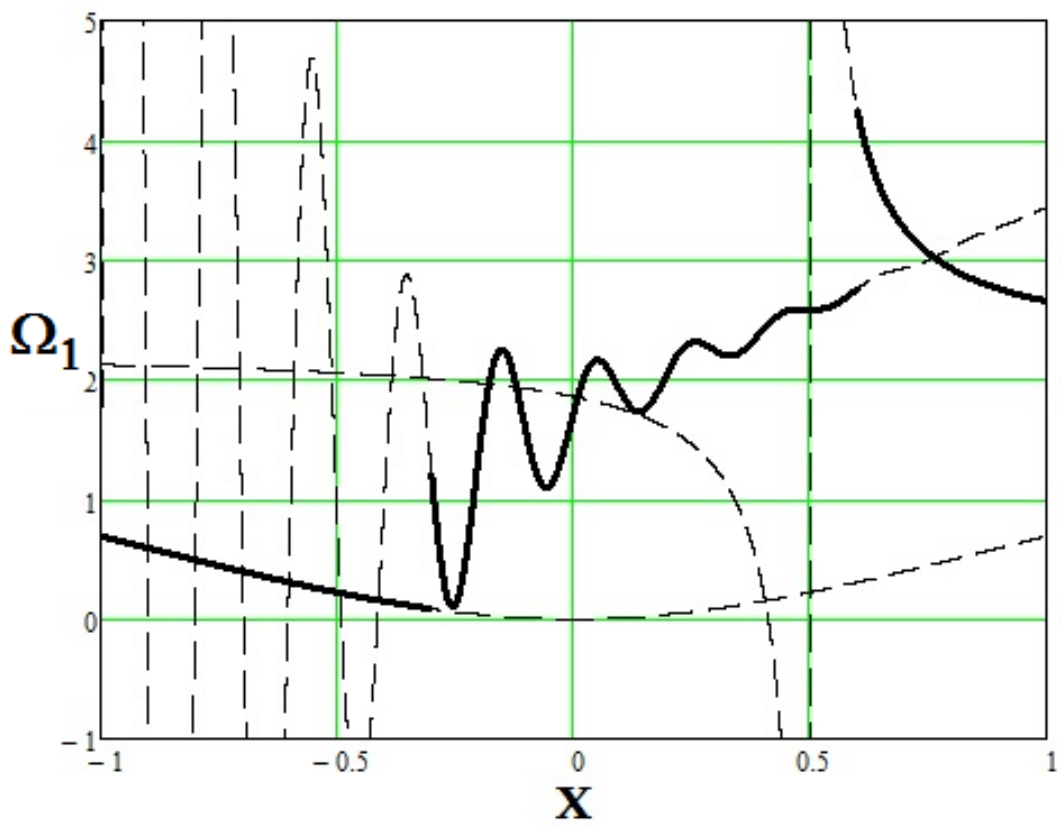

**Fig. 9** The three continuous traces represents our approximating function $\boldsymbol{\Omega_1}$ to the discontinuous function $\boldsymbol{\Psi_1}$ plotted in Fig. 2. The dashed curves represent the three partition functions $\Psi_j$ , *j=1, 2, 3,* that define $\boldsymbol{\Psi_1}$. Compare it with Fig. 10. Note the accuracy of our approximation, particularly at the two discontinuities, and the smooth joints of the three approximant traces with the plots of the three separate partition functions. Also note that our approximant overcomes the divergence of the hyperbolic partition function $\Psi_3$ at $x=0.5$.

The three continuous traces in Fig. 9 represent our approximating function $\boldsymbol{\Omega_1}$ to the discontinuous function $\boldsymbol{\Psi_1}$. The dashed curves are true plots of three partition functions that define that function in the entire domain of the given discontinuous function. The smooth joints of our approximant traces with the three original partition functions at the discontinuity points is something to be noted, it is indeed a sign of the good quality of our approximant. It may also be noticed that our approximant also overcomes the divergence at $x=0.5$ of the third partition function (a hyperbola) which presents a vertical asymptote at that point. Thus the plot of our approximant does faithfully represent the given discontinuous function.

In Fig.10 we have plotted our approximant $\boldsymbol{\Omega_1}$ in a truly small neighbourhood of size $2\times10^{-10}$ centred at the discontinuity point $x=-0.3$ of $\boldsymbol{\Psi_1}$**.** It may be seen that our approximating function faithfully reproduces the jump of that function (see Fig. 8) at that point without defects such as *overshooting* or the Gibbs phenomenon. In this sense our method of approximation is better than any of the methods of approximation based on truncated series or Padé polynomials. In passing, note that the neighbourhood about $x_1 = -0.3$ depicted in this figure is so small that the approximant $\boldsymbol{\Omega_1}$ when plotted there necessarily appear as horizontal traces (no as curved traces about $x=-0.3$ as in Fig. 8).

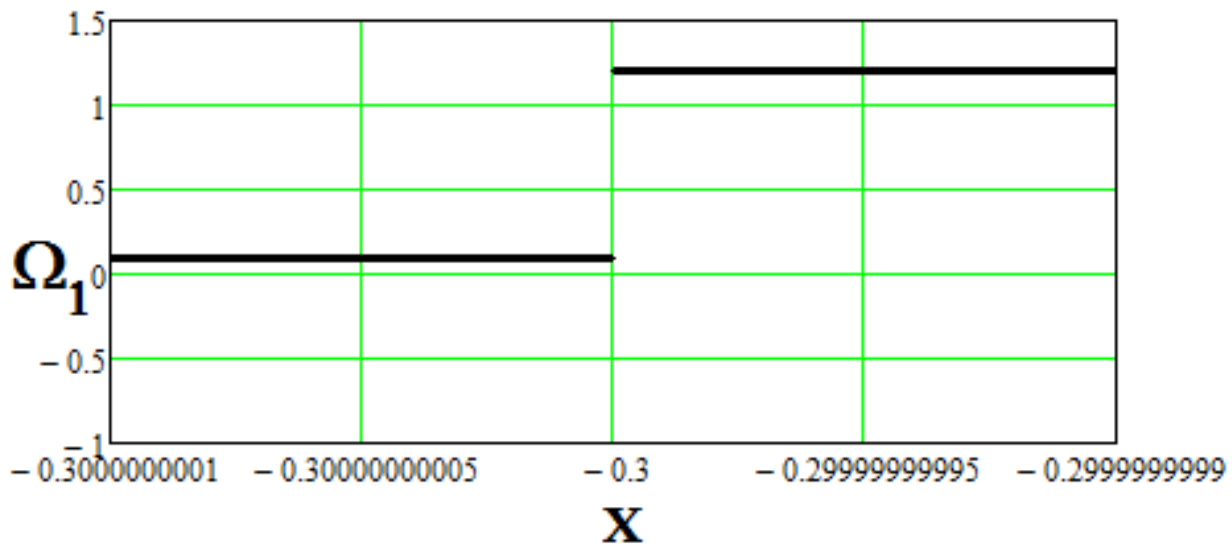

**Fig.10** Plot of our approximant $\boldsymbol{\Omega_1}$ in a truly small neighbourhood of width $2\times10^{-10}$ about the discontinuity point x= - 0.3 of the discontinuous function $\boldsymbol{\Psi_1}$**.** The "jump" there may be seen in Fig. 8, and it is faithfully reproduced by our approximant without defects whatsoever.

In Figs. 11 (a), (b) and (c) we have plotted the relative error $\boldsymbol{e(x)} = (\Omega_1 - \boldsymbol{\Psi_1})/\boldsymbol{\Psi_1}$ of the approximant obtained by applying our method to the current example in the three domains of the partition functions of the discontinuous function Ψ1: (-1,-0.3), (-0.3, 0.6) and (0.6, 1), respectively. It may be seen that the relative errors are truly negligible.

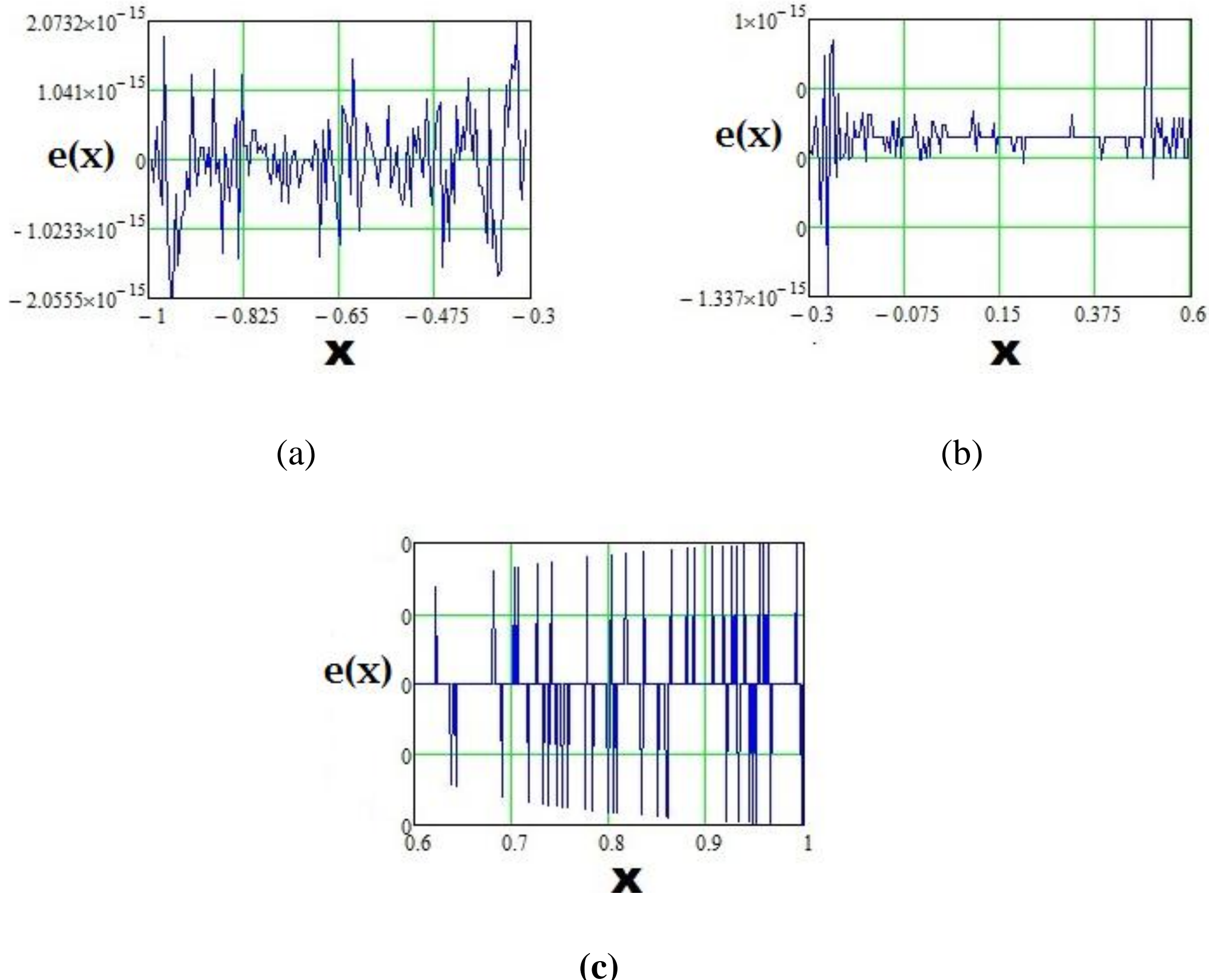

(a) (b)

**(c)**

**Fig. 11** The relative error $e(x) = (\Phi_1 - \Psi_1)/\Psi_1$ of our approximant $\boldsymbol{\Phi_1}$ to the discontinuous function $\boldsymbol{\Psi_1}$ in the partition domains of the latter: (a) in (-1, -0.3) (b) in (-0.3, 0.6), and (c) in (0.6, 1). It may be seen that the maximum relative error is less than $1.1\times10^{-14}$ except at the point $x= 0.5$ where the third partition function presents a vertical asymptote (however as shown in Fig. 12 the relative error there is still negligible)

We considered to be interesting to show the value of the relative error of our approximant at the point $x=0.5$ in the domain (-0.3, 0.6) of the third partition function $\Psi_3$. It may be recalled that the latter represents a hyperbola with a vertical asymptote of equation $x=0.5$. Fig. 11 (b) shows that such local relative error is much larger than the relative error elsewhere in that domain (errors of the order of $10^{-15}$). As shown in Fig. 12 the pathological singularity at $x=0.5$ produces in fact a local relative error of only about 0.3 % of our approximant, an error which of course may be considered negligible for most purposes.

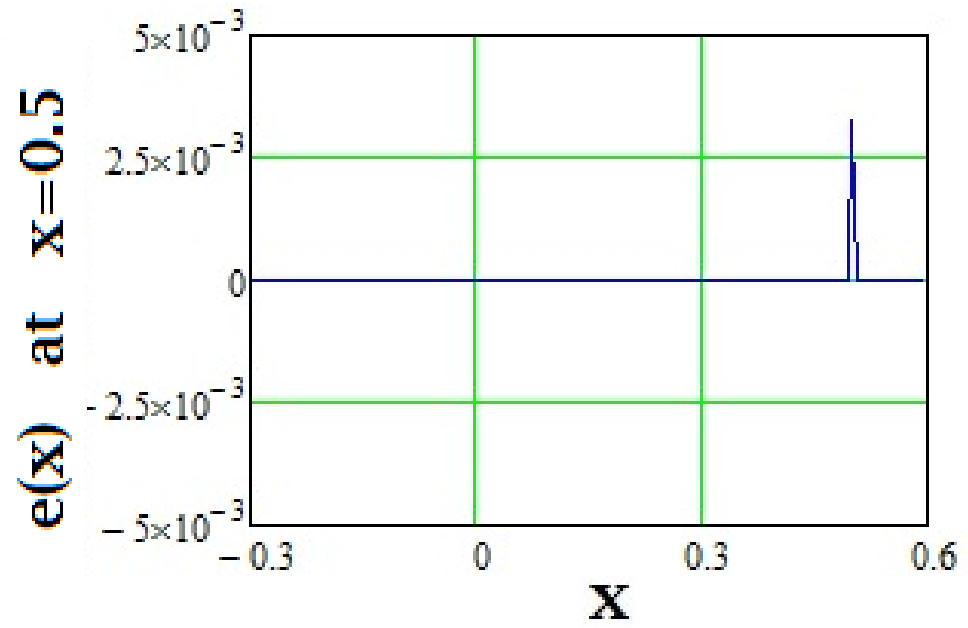

**Fig. 12** Local relative error of the approximant at $x=0.5$. The partition function in the interval (-0.3, 0.6) is a hyperbola whose asymptote is the vertical line $x=0.5$. In spite of this singularity $x=0.5$ the local relative error of the approximation is negligible, only about 0.3 %.

## Example IV

In this example we introduce a discontinuous thin *rect* function $\boldsymbol{\Theta}$ of small width $w=2h$ ($h=1\times10^{-5}$), whose domain is ($x_0=-1$, $x_f=1$) and of height 1/(2h), which when plotted looks like a Dirac $\delta$-function. This function is centred at the point $q=0.3333$, and therefore discontinuous at the points $x_1 = q - h$ and $x_2 = q + h$. It shall be seen that our approximant $\boldsymbol{\Omega}$ to this function does behave as a *Dirac* $\delta$-like function. But, let us first obtain the approximant to our thin *rect* function $\boldsymbol{\Theta}$.

The three partition functions of the given discontinuous function $\boldsymbol{\Theta}$ and their respective domains are of course:

$$\sigma_1(x)=0 \text{ if } x\in(x_0, x_1); \quad \sigma_2(x)=1/(2h) \text{ if } x\in[x_1, x_2]; \quad \sigma_3(x)=0 \text{ if } x\in(x_2, x_f). \quad (34)$$

We now resort to Eq. 1 to define the two connecting functions $\chi_1$, $\chi_2$ at the discontinuity points $x_1$ and $x_2$. With these connecting functions at hand we proceeded to calculate the three auxiliary functions $F_1= 0$, $F_2= 1/(4h)$, $F_3= -1/(4h)$, using the procedures presented in Section 2. Thus our approximant to the given *rect* function is:

$$\Omega(x)=1/(4h)[\chi_1(x)- \chi_2(x)],$$

which appears plotted below in Fig. (13), and does look like a *Dirac delta* function.

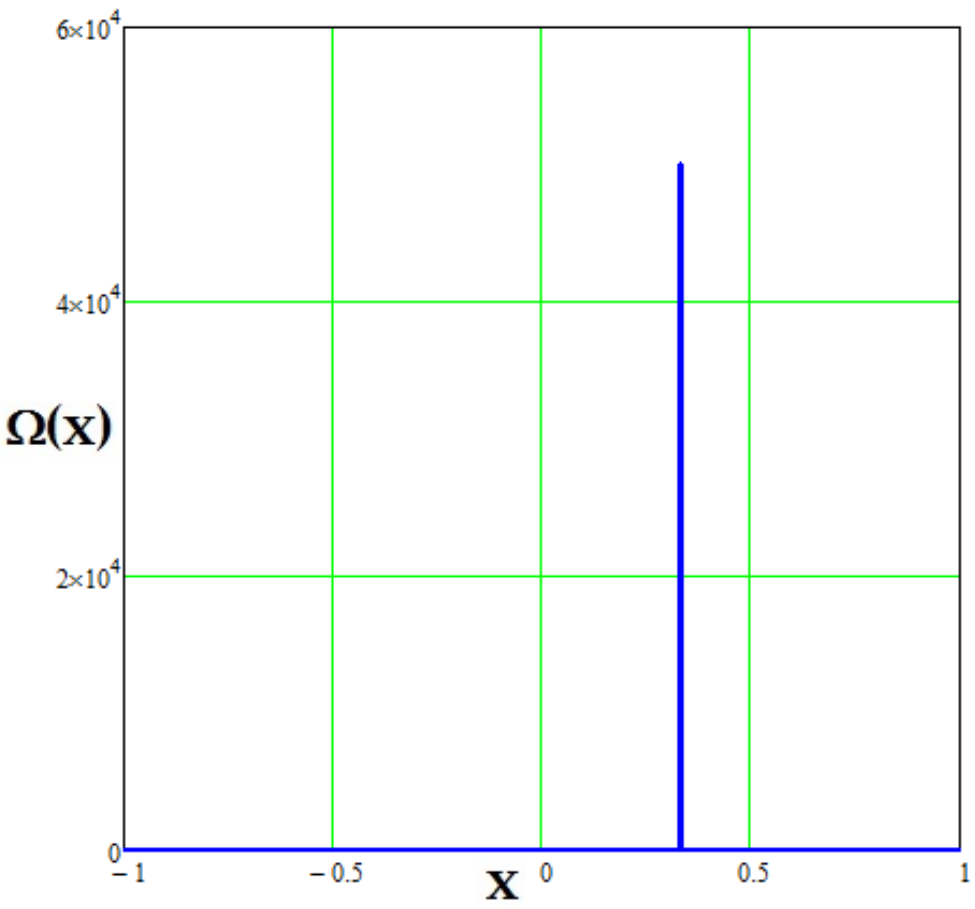

**Fig. 13** Approximant to the thin rect function of width $w=0.00004$ with discontinuities at $x_1=0.33332$ and $x_2= 0.33334$ and defined by Eqs. (34) As shown in the text this function behaves as a *Dirac* $\delta$*-like* function.

Let us now show that this approximant behaves as if it were a *Dirac* $\delta$-like function. In effect:

(i) By definition a *Dirac*-$\delta$ function satisfies the relation: $\int_{x_0}^{x_f} \delta(x-q)dx = 1$. It is therefore meaningful to compute this integral for our approximant $\Omega$. The result is:

$$I_1 = \int_{x_0}^{x_1} \Omega(x)dx + \int_{x_1}^{x_2} \Omega(x)dx + \int_{x_2}^{x_f} \Omega(x)dx = 1.0000000000$$

which is equal to the expected value 1 within the truly small (computed) relative error $e_1$:

$$e_1 = (I_1 - 1)/1 = 9.999\times 10^{-13}.$$

(ii) Let us now consider the important *Sifting Property* of a Dirac $\delta$-function, namely: $\int_{x_0}^{x_f} \delta(x-q)f(x)dx = f(q)$ for any function *f*. Moreover, let us consider the particular case of the sine function and write an analogous sifting relation for our approximant $\Omega$:

$$I_2 = \int_{x_0}^{x_f} \Omega(x-q)\sin(x)\,dx$$

A direct computation of this integral gives (recall $q= 0.33333$ is given above):

$$I_2 = \int_{x_0}^{x_1} \Omega(x-q)\sin(x)\,dx + \int_{x_1}^{x_2} \Omega(x-q)\sin(x)\,dx + \int_{x_2}^{x_f} \Omega(x-q)\sin(x)\,dx = 0.3271946968\,,$$

while the expected value of this integral when the actual Dirac $\delta$-function is used is of course $\sin(q) = 0.3271946968$, very close to the value obtained with our approximant. It is therefore pertinent to calculate the relative error $e_2$ when its approximant $\Omega$ is used:

$$e_2 = (I_2 - \sin(q))/\sin(q) = -1.567\times 10^{-11}.$$

The two small relative errors $e_1$, $e_2$ computed above show that our approximant is a truly accurate one, showing it is a good replacement for a true Dirac $\delta$-function.

**Example V An approximant for a discontinuous force**

As already mentioned in the Introduction the dynamics of a body subjected to a discontinuous force that depends either on position or time is a relevant problem in mathematical-physics. It leads to ordinary motion differential equations with discontinuous independent terms [14, 16, 20]. Here we consider an example of such discontinuous dynamics: a 1-dimensional mass-spring $(m, k)$ oscillator perturbed by an external time-dependent force that is equal to a constant $f_0$, during a finite time interval $(t_1, t_2)$ comparable to the period $T$ of the free oscillations of the system, and nil otherwise, i.e.

$$f(t) = f_0 \quad \text{if} \quad t \in [t_1, t_2]\,, \tag{25a}$$

$$f(t) = 0 \quad \text{if} \quad t \notin [t_1, t_2], \tag{25b}$$

a problem initially posed and solved by Landau and Lifschitz using an alternative special method [20].

The general motion equation of this forced harmonic oscillator is of course

$$m\,\ddot{x}(t) + k\,x(t) = f(t). \tag{26}$$

This is in fact a good example of the important class of ordinary differential equations with discontinuous right-hand sides studied by Filippov [14]. Here we intend to apply our approximation method to find first an approximant $\boldsymbol{\Omega(t)}$ for the discontinuous force $f(t)$ acting on the oscillator, and then to solve its motion equation yet after replacing the force in the righthand sides of Eq. (26) by the approximant.

Consider then the following particular case: the oscillations of a frictionless mass-spring oscillator of mass $m=1$ $[kg]$ and a spring constant $k=4$ $[N/m]$ are perturbed by a constant force of magnitude $f_0$= *4 [N]* during the time interval from $t_1$=*2.00 [s]* to $t_2$=*16.65* [s], being zero otherwise. If we take the arbitrary time domain from $t_0$ = *-10 [s]* to $t_f$ = *30 [s]*, the three partition intervals of the discontinuous force are of course

*(-10, 2)*, *[2, 16.65]*, and *(16.65, 30)*. According to Eqs. (25), the three partition functions of the applied discontinuous force are therefore the constant forces *f* = *0 [N]*, *f*= *4 [N]* and *f*=*0 [N]*, respectively.

Once the two connecting functions at $t_1$ = *2 [s]* and $t_2$ =*16.65 [s]* are defined using Eq. (1) the nine required $\{S_{jn}\}$ coefficients are once again readily calculated using Eqs. (7) of Section 2: $\{S_{00}=1, S_{10}=1, S_{20}=1, S_{01}=-1, S_{02}=-1, S_{11}=1, S_{12}=-1, S_{21}=1, S_{22}=1\}$. Then, using Eq. (9) the three auxiliary functions of the sought approximant are found to be: $F_0$ = *0 [N]*, $F_1$ = *2 [N]* and $F_2$ = *-2 [N]*. With these at hand we may write our approximant $\boldsymbol{\Omega}$ to the discontinuous force acting on the mass-spring system as

$$\boldsymbol{\Omega}(t) = F_1 \, tanh\left[\left[\left(\frac{t-t_0}{t_f-t_0}\right)^2\left(\frac{t-t_1}{t_f-t_0}\right)\left(\frac{t_f-t}{t_f-t_0}\right)^2\right]^{-1}\right] - F_2 \, tanh\left[\left[\left(\frac{t-t_0}{t_f-t_0}\right)^2\left(\frac{t-t_2}{t_f-t_0}\right)\left(\frac{t_f-t}{t_f-t_0}\right)^2\right]^{-1}\right] \qquad (28)$$

or

$$\boldsymbol{\Omega}(t) = 2 \, tanh\left[\left[\left(\frac{t+10}{40}\right)^2\left(\frac{t-2}{40}\right)\left(\frac{30-t}{40}\right)^2\right]^{-1}\right] - 2 \, tanh\left[\left[\left(\frac{t+10}{40}\right)^2\left(\frac{t-16.65}{40}\right)\left(\frac{30.-t}{40}\right)^2\right]^{-1}\right] \qquad (29)$$

that we have plotted in Fig. 14. Note that as in the previous examples, the plotted approximant again looks practically indistinguishable from the true force applied to the oscillator given by Eqs. (25).

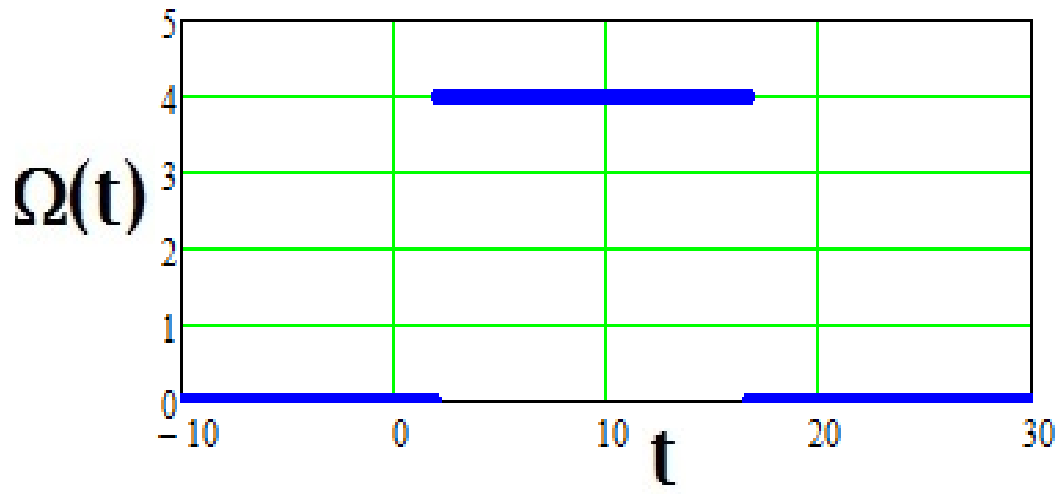

**Fig. 14** Plot of our approximant ***Ω(t)*** to the discontinuous force *f(t)* applied to the oscillator example ($f_0$=*4[N]*)

As announced above this approximating function $\boldsymbol{\Omega}(t)$ we now introduce into the righthand side of the motion differential equation of our mass-spring system (for which $m=1$ *kg*, $k= 4$ *N/m*) to obtain an approximated form of the oscillator motion equation:

$$m\,\ddot{x}(t) + k\,x(t) = \boldsymbol{\Omega}(t). \qquad (30)$$

This is a new differential equation, similar to Eq. (27), that we can solve in a number of ways. For instance, we have obtained the exact analytical solution of this differential equation and found the solution $x(t)$ to be given in terms of our approximant $\boldsymbol{\Omega}(t)$ by:

$$x(t) = c_1 \cos(qt) + c_2 \sin(qt) + \frac{\sin(qt)}{q}\int_{t_0}^{t} \boldsymbol{\Omega}(t')\cos(qt')dt' + (-1)\frac{\cos(qt)}{q}\int_{t_0}^{t} \boldsymbol{\Omega}(t')\sin(qt')\,dt,$$

where $q=2$, $c_1=x_0$, $c_2= v_0/q$. This solution we have plotted in Fig. 15 below, for the initial conditions $x_0 =0$, $v_0=0$.

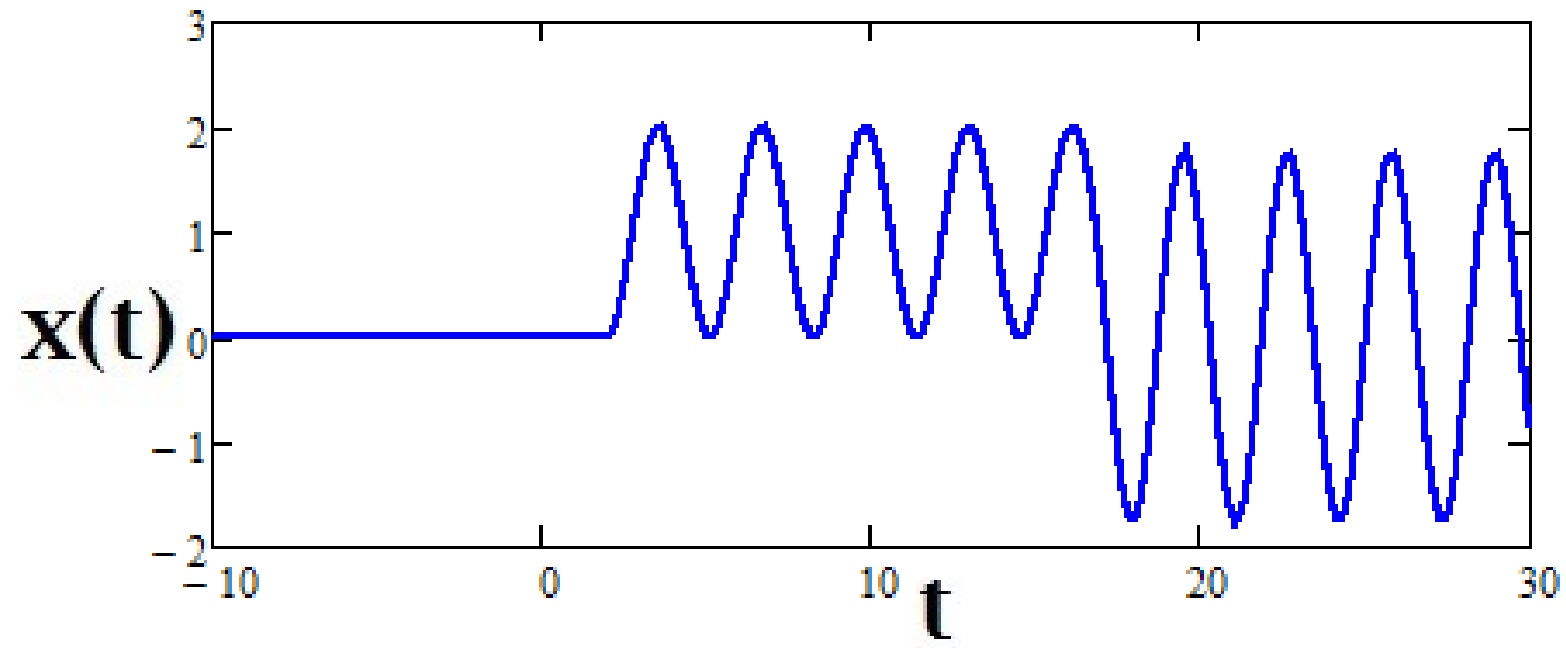

**Fig. 15** Oscillations of the point mass about $x=1$ as predicted by a motion differential equation where the actual impulsive force function $f(t)$ was replaced by our approximant $\boldsymbol{\Omega}(t)$ in Eq. (29). Note that as the applied force goes to zero at $t_2 = 16.65$ the oscillations take place about the expected new equilibrium position $x=0$.

The oscillating function plotted in Fig. 15 is indistinguishable from the true oscillations of our mass-spring system that we also obtained by solving the exact motion Eq. (26) of the mass-spring system. It may be noticed that since the force interaction $f(t)$ is only applied during the finite interval $[2.00, 16.65]$ the amplitude of oscillation in that interval is smaller than the amplitude of the free oscillations plotted at the right of of Fig. 15. Also note that during the force interaction, in the interval $[2, 16.65]$, the point mass oscillates about $x=1$ with amplitude $1$, and then as the interaction ceases at $t= 16.65$ it freely oscillates about $x=0$ with amplitude $1.8$, exactly as can be obtained using the special method of Landau and Lifschitz [20]. Actually, the oscillation amplitude for $t > t_2$ depends upon the exact instant $t_2$ in which the external force $f$ goes to zero.

A similar case of forced oscillation was recently presented by Donoso and Ladera [24] in which the driven vertical oscillations of a conducting ring suspended from a spring are suddenly perturbed by an oscillating magnetic field of high frequency (produced by a surrounding coil that is excited with an AC current during a finite time) that results in a constant average vertical magnetic force on the ring. It was then predicted and experimentally measured (Fig. 16) that the ring equilibrium position and its amplitude of oscillation behaved as the function plotted in Fig. 15.

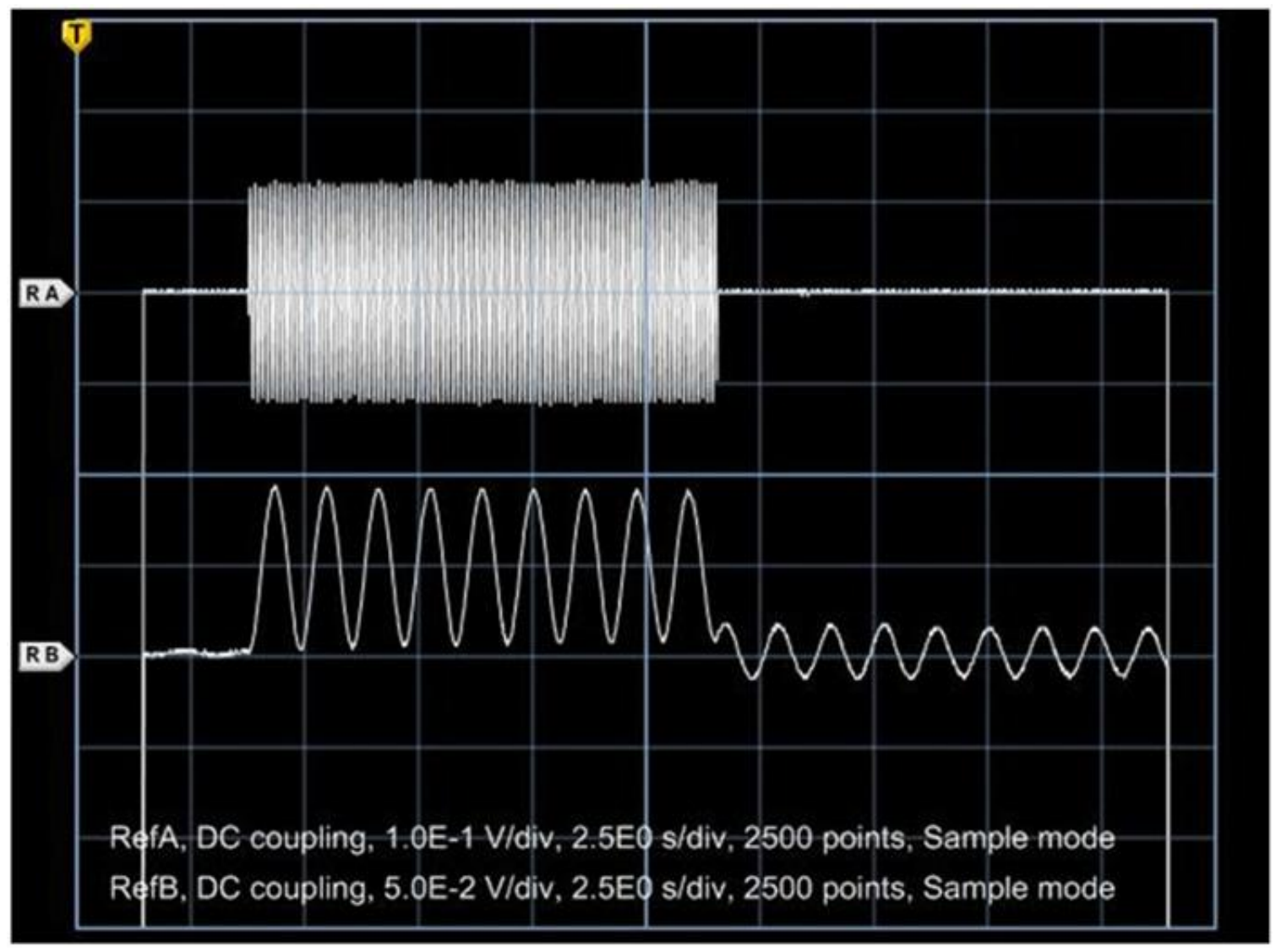

**Fig. 16** Measured vertical oscillations of a conducting ring suspended from a weak spring (lower oscilloscope trace) driven by an AC field. When an additional constant magnetic force displaces upward the ring the equilibrium position is raised (larger amplitude oscillations in the left of the trace), as soon as that force ceases the ring returns to its natural oscillations about its initial equilibrium position (lower amplitude oscillations in the right of the trace) ([24] doi: 10.1088/0143-0807/35/1/015002, *Reproduced by permission of IOP Publishing. All rights reserved*)

## IV Discussion and Conclusions

A new and very accurate method to approximate piecewise continuous functions, with a finite number of discontinuities and partition functions, has been presented in this work. The method was then applied to a set of simple yet interesting mathematical examples for the reader to grasp the essentials of the method and how to apply it. Our method relies on the definition of *connecting functions* at each discontinuity point, functions that are defined in terms of hyperbolic tangent functions whose arguments are rational functions of Newton binomials. As described in details in Section 2, our method consists essentially in creating an approximating function in terms of those connecting functions, actually a linear combination of such functions, in which the coefficients are functionals found by solving a linear system of inhomogeneous equations, and where the role of independent terms is taken by the partition functions of the initially given piecewise continuous function. We then applied our method to two well-known discontinuous functions, and to a function with two discontinuities and whose three partition functions were rather arbitrarily defined (the third of them having a vertical asymptote in its domain for us to test whether or not our approximant were to fail at the location of that asymptote). We also illustrated our method considering the simple yet interesting case of a mass-spring oscillator that is forced by a constant force during a finite period of time; we then used our method to find an approximant that accurately represented the actual oscillations of that system. This successful application of our method to a forced oscillator is relevant because it opens a door to its application to physics systems represented by ordinary differential equations whose right-hand side is a discontinuous function. We showed that our method generates approximants to discontinuous functions that are very accurate, giving relative errors of the order of $10^{-12}$ and even less than that, and that they do converge at the discontinuity points much better that approximants obtained using well-known alternative methods based on truncated series or *e.g.*Padé-Jacobi polynomials. Our method can therefore be applied to interesting cases in physics such as the so-called Filippov´s problems, Fresnel diffraction, propagation of waves in layered dielectrics, and to the WKB method, as well as to statistical physics and physical chemistry problems such as the Van der Waals isotherms and discontinuities in first order thermodynamics transitions. We have also obtained an approximated solution for the important nonlinear Riccati´s differential equation $dy/dx + y^2 = R(x)$ when the right-hand function $R$ is a piecewise continuous function with two discontinuities, and generalized the approximation to any number of discontinuities of such function. This particular case of Riccati differential equation and an additional number of similar discontinuous function problems we have already solved using our method of approximation and will be published elsewhere. Finally, in Appendix A we showed how to regularize our connecting functions by introducing a small parameter $\sigma^2$ in the denominator of the argument of such functions. With that simple redefinition we warrant our approximants (to discontinuous functions) to be continuous everywhere and regular enough to be successfully replaced into Filippov´s differential equations.

**Appendix A: Regularizing the connecting functions**

As described in Section II our method of approximation rests on the elementary definition (Eq. (1)) of the connecting functions $\chi$ between neighbouring partitions functions. As there explained the method amounts to a rather simple, short and precise procedure that leads to a representation of any discontinuous function in terms of such elementary functions. Yet, when closely examined one discovers that such connecting functions are not regular enough at the discontinuity points themselves: at such points the argument of the hyperbolic tangent in Eq. (1) diverges and the connecting functions take the finite values ±1, i.e they become discontinuous. In this Appendix we redefine the form of our connecting function $\chi$ in Eq. (1) to solve this fundamental difficulty. We shall denote the new connecting function as $\Gamma$ and be referred to as the *regularized connecting function*, whose most relevant features will be their continuity at the discontinuity points of the approximated discontinuous function, and that its recovery of the connecting functions as given by Eq. (1) under the proper conditions. Thus the improved regularized connecting function will allow us fully recovering the discontinuity jumps, yet using a continuous function (see Fig. A-1). Also note that our *regularized* connecting function will be particularly useful when solving the motion equations of discontinuous dynamical systems, the so-called Filippov´s problems.

Our new *regularized connecting function* $\Gamma$ at each discontinuity point $x_d$ we have chosen to define as:

$$\Gamma(x;x_d,\sigma) \equiv \tanh\left[\frac{(x_f-x_0)^5(x-x_d)}{[(x-x_0)(x-x_d)(x-x_f)]^2+(x_f-x_0)^4\sigma^2}\right], \qquad \text{(A-1)}$$

where $\sigma$ is a small real number, that as explained below is related to the small width of the jump of the connecting function at the discontinuity point. The role of this parameter $\sigma$ is rather important: a small value of it, say $\sigma \in (10^{-4}, 10^{-20})$, means that one will get sudden jumps of the connecting function $\Gamma$ at all discontinuity points, while a reasonably larger value of it, say $\sigma \geq 10^{-3}$ will give a smoother continuous jump in a rather small interval about each discontinuity point, *i.e.* these latter values of the parameter $\sigma$ might become useful if one needs to joint two neighbouring partitions functions in a smooth way. Of course most of the time one needs this free parameter to be rather small say $\sigma \leqq 1.10^{-4}$. Moreover, and what is truly important: *our regularized connecting function is always well defined and shows no Gibbs ringing phenomenae at all* (see Fig. (A-1)). It is straightforward to check that the initial connecting function $\chi$ of Eq.(1) can be recovered from our regularized connecting function $\Gamma$ of Eq. (A-1) by simply letting $\sigma$ to be negligible (say $\sigma \in (10^{-4}, 10^{-20})$).

The small interval where the actual and appropriate behaviour or profile of our regularized connecting function can be observed is about $(x_d - \sigma^2/L,\ x_d + \sigma^2/L)$ where $L = x_f - x_0$ is the length of the discontinuous function domain. Fortunately, $\sigma^2$ is always the size of a truly small neighbourhood about the abscissa $x_d$. This we can show in the *zoomed* and rather meaningful plots in Fig A-1. There we have considered (just for the sake of illustration) a certain discontinuous function defined in the domain $(x_0=0,\ x_f=2)$ with a single discontinuity point at $x_d=0.3$, and have plotted the two regularized connecting functions $\Gamma(x;\sigma)$ in the interval $(x_d - 1\times10^{-4},\ x_d + 1.0\times10^{-4})$ about $x_d$, for two different values of the parameter, namely $\sigma=5.0\times10^{-3}$ and $\sigma=5.0\times10^{-6}$ in Eq. (A-1). While the figure shows a smooth rising dotted curve for the larger parameter value $\sigma=5.0\times10^{-3}$ in the small interval of length $\sigma^2/L=2.5\times10^{-5}$, it also shows a neat square jump (solid line) for the smaller parameter $\sigma$ in the same interval. The comparison is straightforward and the difference is indeed significant.

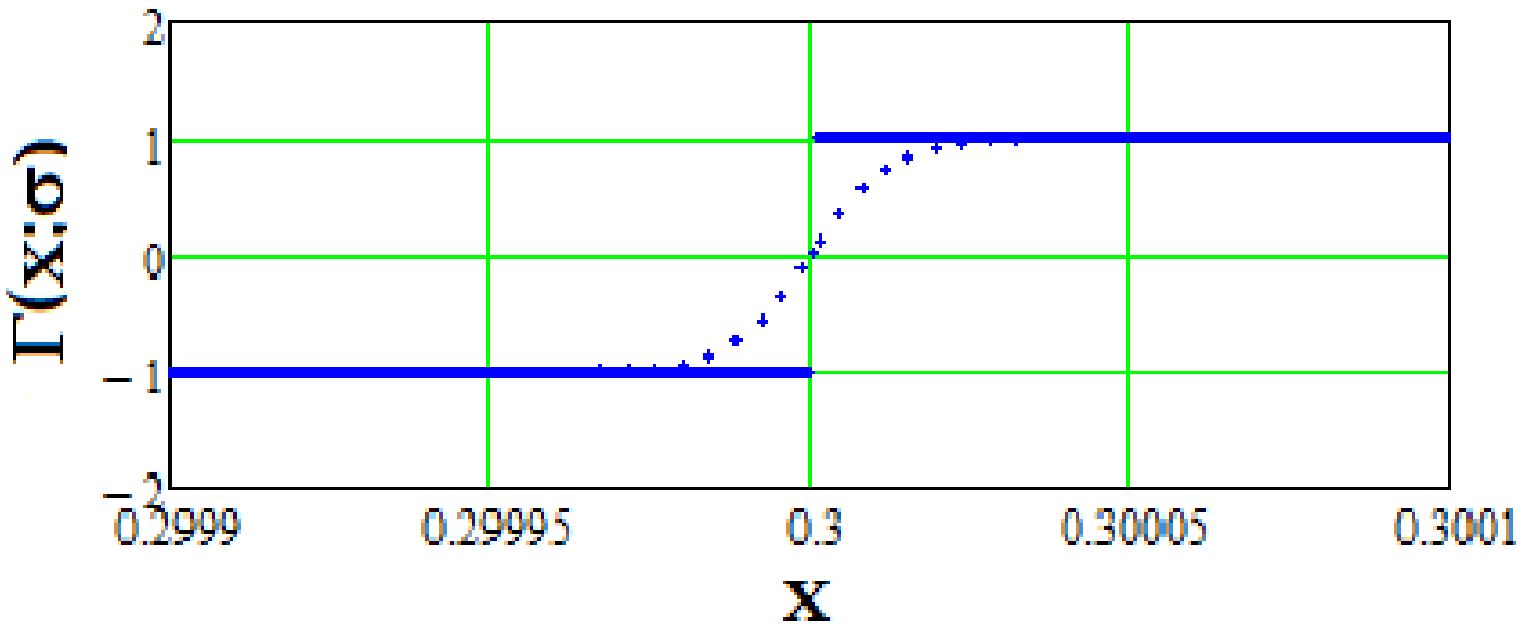

**Fig. A-1** Plots of the regularized connecting functions $\Gamma(x;\sigma)$ in the small interval $(x_d - 1.0\times10^{-4},\ x_d + 1.0\times10^{-4})$ about the discontinuity point at $x_d=0.3$ for $\sigma=5.0\times10^{-3}$ (dotted line) and for $\sigma=5.0\times10^{-6}$ (solid line). As expected the first connecting function shows a smooth rising behaviour; the second connecting function shows a sudden jump.

## Appendix 2

In this appendix we present the derivative $\Gamma_{D,sld}$ of our connecting function in Eq.(1) of the main text:

$$\Gamma_{D,sld} = \mathrm{sech}\left[\frac{(x_f - x_i)^5 (x - x_0)}{f(x)}\right]^2 (x_f - x_i)^5 \left[\frac{f(x) - (x - x_0) Df(x)}{[f(x)]^2}\right]$$

where the auxiliary function *f(x)* in the denominator is the following 6th degree polynomial:

$$f(x) = f_0 + f_1.x + f_2.x^2 + f_3.x^3 + f_4.x^4 + f_5.x^5 + f_6.x^6$$

whose six coefficients $f_i$ are as follows:

$$f_0 = b^2x_0^2x_f^2 + b^2x_i^2x_f^2 + b^2x_i^2x_0^2 + x_i^2x_0^2x_f^2 + b^6 + b^4x_0^2 + b^4x_i^2 + b^4x_f^2,$$

$$f_1 = -2b^4x_f - 2b^4x_i - 2b^2x_ix_0^2 - 2b^2x_ix_f^2 - 2x_i^2x_0x_f^2 - 2b^2x_i^2x_0 - 2x_i^2x_0^2\,x_f - 2b^2x_fx_i^2 - 2b^2x_0x_f^2 - 2b^2x_fx_0^2 - 2x_ix_0^2x_f^2,$$

$$f_2 = -2b^2x_i^2 + 2b^2x_f^2 + 2b^2x_0^2 + x_i^2x_0^2 + x_f^2x_0^2 + 4x_i^2x_0x_f + 4x_f^2x_0x_i + 4b^2x_ix_0 + 4x_ix_0^2x_f + 4b^2x_ix_f + 4b^2x_0x_f + 3b^4 + x_i^2x_f^2,$$

$$f_3 = -2x_i^2x_f - 2x_i^2x_0 - 2x_f^2u_0 - -2x_f^2x_i - 4b^2x_i - 2x_0^2x_f - 4b^2x_0 - 4b^2x_f - 2x_ix_0^2 - 8x_ix_fx_0,$$

$$f_4 = 4x_0x_f + 4x_0x_i + 4x_ix_f + x_i^2 + x_f^2 + 3b^2 + x_0^2,$$

$$f_5 = -2x_0 - 2x_f - 2x_i,$$

$$f_6 = 1\ .$$

This first derivative $\Gamma_{D,sld}$ is plotted in Fig A2. As shown the discontinuity of the function $\chi_{sld}$, in (Eq. (1) was set at $x_0=3$.

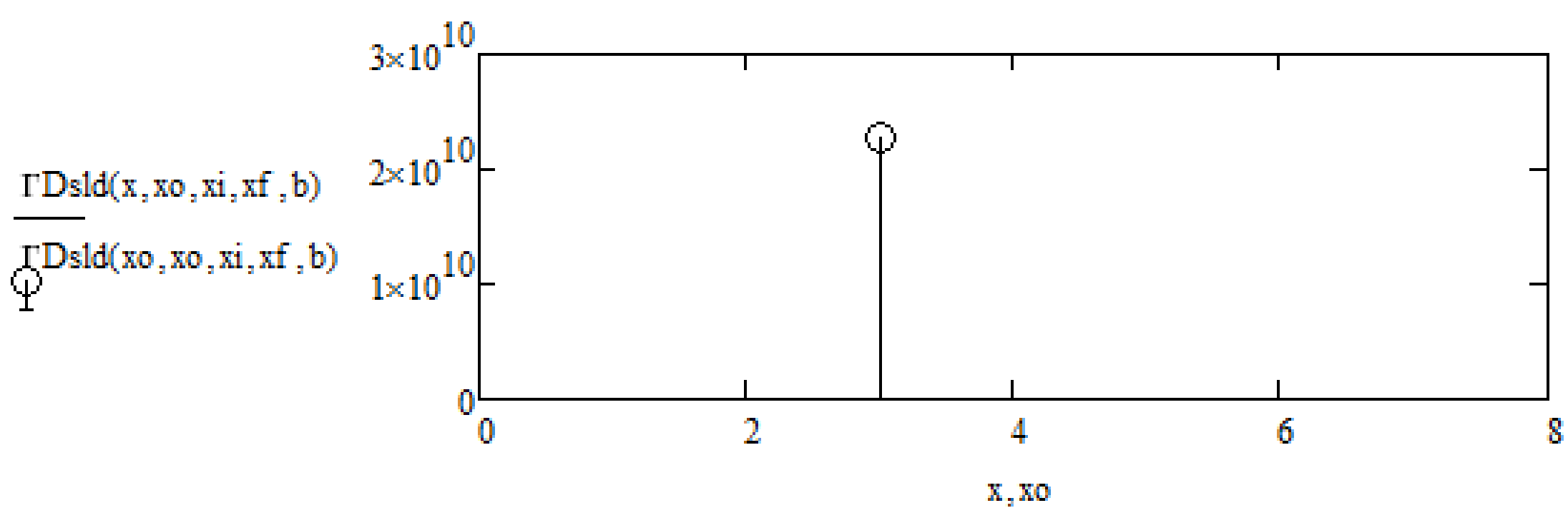

**Fig A2** Plot of the first derivative of our connecting function $\chi_{sld}$ (in Eq. (1)). The single plotted *circle* at the point (3, 2.28×10^10) is the actual high value of the derivative $\Gamma_{D,sld}$ at $x_0=3$

## References

**[1]** K. S. Eckoff Accurate and efficient reconstruction of discontinuous functions from truncated series expansions, *Math. Comp*. **61** (1993) 745-763

**[2]** L. Emmel, S. M. Kaber, Y. Maday, Pade-Jacobi Filtering for spectral approximations of discontinuous solutions, *Numer. Algo*. **33** (2003) 251-264

**[3]** T. Chantrasmi, A. Doostan, G. Iaccarino, Padé–Legendre approximants for uncertainty analysis with discontinuous response surfaces, *J. Comp. Phys.* **228** (2009) 7159–7180

**[4]** J. S. Hesthaven, S.M. Kaber, and L. Lurati, Pade-Legendre Interpolants for Gibbs Reconstruction, *J. Sci. Comp.* **28** (2006) 337-359

[5] D. Costarelli, Ph. D. Thesis, Sigmoidal Functions Approximation and Applications, Universitat degli Study Roma Tres, Roma, Italy (2013)

**[6]** S. M. Kaber and Y. Maday, Pade-Chebyshev approximants, *SIAM Journal on Numerical Analysis*, **43** (2004) 437-454

**[7]** G. Kvernadze, Approximating the jump discontinuities of a function by its Fourier-Jacobi coefficients, *Math. Comp*. **73** (2003) 731-751

**[8]** J. Awrejcewicza, Michal Feˇckan, P. Olejnika On continuous approximation of discontinuous systems, *Nonlinear Analysis* **62** (2005) 1317–1331

**[9]** D. Mumford, J. Shah, Optimal Approximations by Piecewise Smooth Functions and Associated Variational Problems, *Communications on Pure and Applied Mathematics*, Wiley, New York, 1989, Vol. XLII 577-685

**[10]** G. Donoso and P. Martín, Ion-Acoustic Dispersion Relation with Direct Fractional Approximations for Z'(x)", *J. Math. Phys.* **26** (1986) 1186-1188

**[11]** C. L. Ladera and P. Martín Transmittance of a Circular Aperture by an Integrable Fractional-Like Approximation to $J_0(x)$, *J. Comput. Phys.* **73** (1987) 481-489

**[12]** E. Stella and P. Martin, Non-Linear Electromagnetic Field Diffusion for Plasmas with Non-Linear Conductivity by Two Point Quasifractional Approximants, *Astrophys Space Sci.,* 256 (1998) 283-287

**[13]** E. Hecht, *Optics*, second ed., Addison, Reading, Ma.,1990, Ch. 10, pp. 456-458

**[14]** A. M. Filippov, *Differential equations with discontinuous righthand sides,* Kluwer, Dordrecht, 1988

**[15]** M-F Danca and S. Codreanu, On a possible approximation of discontinuous dynamical systems, *Chaos Solitons and Fractals* **13** (2002) 681-691

**[16]** M-F Danca, On a class of discontinuous dynamical systems, *Mathematical Notes, Miskolck* **2** (2001) 103-116

**[17]** R. I. Leine, D. H. Van Campen and D. L. Van de Vrande, Bifurcations in Nonlinear Discontinuous systems, *Nonlinear Dynamics* **23** (2002)105-164

**[18]** K. Huang, Statistical Mechanics, Wiley, New York. 1987, second edition,

**[19]** H. E. Stanley, *Introduction to Phase Transitions and Critical Phenomena*, Osford Univ. Pr. Oxford, 1987.

**[20]** L. D. Landau and E. M. Lifschitz *Course on Theoretical Physics*, Vol. 1: Mechanics, Pergamon, Oxford, 1969, Ch. V pp. 61-65

**[21]** F. M. Ceraglioli, *Discontinuous ordinary differential equations and stabilization*, Tesi di dottorato XI Ciclo, Universita Firenze, Italy (1999)

**[22]** K. al Shanmari, *Filippov´s Operators and Discontinuous Differential Equations* Ph. D. Thesis, Louisiana St. University, USA (2006)

**[23]** C. Kunush, *Control de modos deslizantes de segundo orden de pilas de combustible PEM*, VIII Symposium on Advanced Control of Energy Systems, Almería, Spain (2010)

**[24]** G. Donoso, C. L. Ladera, The naked toy model of a jumping ring, *Eur. J. Phys*. **35** (2014) 015002 doi: 10.1088/0143-0807/35/1/015002